\documentclass[11pt]{amsart}
\usepackage{amsmath}
\usepackage{amsfonts}
\usepackage{amssymb}
\usepackage{amsthm}
\usepackage{mathtools} 
\usepackage{latexsym}

\usepackage{enumerate}
\usepackage{cancel}
\usepackage{cases}
\usepackage{empheq}
\usepackage{multicol}

\usepackage{subfigure}
\usepackage{float}

\usepackage{wrapfig}

\usepackage{color}

\usepackage[color=white,linecolor=black]{todonotes}

\usepackage{mathrsfs}
\usepackage{fontenc} 
\usepackage{inputenc} 

\usepackage{verbatim}

\theoremstyle{plain}
\newtheorem{theorem}{Theorem}[section]

\theoremstyle{definition}

\theoremstyle{remark}
\newtheorem{remark}[theorem]{Remark}

\numberwithin{equation}{section} 
\numberwithin{figure}{section}   

\usepackage[square,comma,numbers,sort&compress]{natbib}
\usepackage{url}

\newcommand{\field}[1]{\mathbb{#1}}

\newcommand{\nZ}{\field{Z}}

\newcommand{\nR}{\field{R}}

\newcommand{\maps}{\rightarrow}

\newcounter{my_counter}
\title[Nonlinear Data Assimilation]{Nonlinear Continuous Data Assimilation}

\date{\today}

\author{Adam Larios}
\address[Adam Larios]{Department of Mathematics, 
                University of Nebraska--Lincoln,
        Lincoln, NE 68588-0130, USA}
\email[Adam Larios]{alarios@unl.edu}
\author{Yuan Pei}
\address[Yuan Pei]{Department of Mathematics, 
                University of Nebraska--Lincoln,
        Lincoln, NE 68588-0130, USA}
\email[Yuan Pei]{ypei4@unl.edu}
\thanks{MSC 2010 Classification: }

\begin{document}
\begin{abstract}
 We introduce three new nonlinear continuous data assimilation algorithms.  These models are compared  with the linear continuous data assimilation algorithm introduced by Azouani, Olson, and Titi (AOT).  As a proof-of-concept for these models, we computationally investigate these algorithms in the context of the 1D Kuramoto-Sivashinsky equation.  We observe that the nonlinear models experience super-exponential convergence in time, and converge to machine precision significantly faster than the linear AOT algorithm in our tests.  
\end{abstract}

\maketitle
\thispagestyle{empty}

\noindent
\section{Introduction}\label{secInt}
\noindent

Recently, a promising new approach to data assimilation was pioneered by Azouani, Olson, and Titi \cite{Azouani_Olson_Titi_2014,Azouani_Titi_2014} (see also \cite{Cao_Kevrekidis_Titi_2001,Hayden_Olson_Titi_2011,Olson_Titi_2003} for early ideas in this direction).  This new approach, which we call AOT data assimilation or the linear AOT algorithm, is based on feedback control at the partial differential equation (PDE) level, described below.  In the present work, we propose several nonlinear data assimilation algorithms based on the AOT algorithm, that exhibit significantly faster convergence in our simulations; indeed, the convergence rate appears to be super-exponential.

Let us describe the general idea of the AOT algorithm.  Consider a dynamical system in the form,
\begin{empheq}[left=\empheqlbrace]{align}\label{ODE}
\begin{split}
\tfrac{d}{dt} u &= F(u),\\
u(0)&=u_0.
\end{split}
\end{empheq}
For example, this could represent a system of partial differential equations modeling fluid flow in the atmosphere or the ocean.  A central difficulty is that, even if one were able to solve the system exactly, the initial data $u_0$ is largely unknown.  For example, in a weather or climate simulation, the initial data may be measured at certain locations by weather stations, but the data at locations in between these stations may be unknown.  Therefore, one might not have access to the complete initial data $u_0$, but only to the observational measurements, which we denote by $I_h(u_0)$. (Here, $I_h$ is assumed to be a linear operator that can be taken, for example, to be an interpolation operator between grid points of maximal spacing $h$, or as an orthogonal projection onto Fourier modes no larger than $k\sim 1/h$.)
Moreover, the data from measurements may be streaming in moment by moment, so in fact, one often has the information $I_h(u)=I_h(u(t))$, for a range of times $t$.  Data assimilation is an approach that eliminates the need for complete initial data and also incorporates incoming data into simulations.  Classical approaches to data assimilation are typically based on the Kalman filter.  See, e.g., \cite{Daley_1993_atmospheric_book,Kalnay_2003_DA_book,Law_Stuart_Zygalakis_2015_book} and the references therein for more information about the Kalman filter.  In 2014, an entirely new approach to data assimilation---the AOT algorithm---was introduced in \cite{Azouani_Olson_Titi_2014,Azouani_Titi_2014}.  This new approach overcomes some of the drawbacks of the Kalman filter approach (see, e.g., \cite{Biswas_Hudson_Larios_Pei_2017} for further discussion).  Moreover, it is implemented directly at the PDE level.  The approach has been the subject of much recent study in various contexts, see, e.g., \cite{Albanez_Nussenzveig_Lopes_Titi_2016,Altaf_Titi_Knio_Zhao_Mc_Cabe_Hoteit_2015,Bessaih_Olson_Titi_2015,Biswas_Martinez_2017,Farhat_Jolly_Titi_2015,Farhat_Lunasin_Titi_2016abridged,Farhat_Lunasin_Titi_2016benard,Farhat_Lunasin_Titi_2016_Charney,Farhat_Lunasin_Titi_2017_Horizontal,Foias_Mondaini_Titi_2016,Gesho_Olson_Titi_2015,Jolly_Martinez_Titi_2017,Markowich_Titi_Trabelsi_2016_Darcy,Mondaini_Titi_2017}.

The following system was proposed and studied in \cite{Azouani_Olson_Titi_2014,Azouani_Titi_2014}:  
\begin{empheq}[left=\empheqlbrace]{align}\label{AOT}
\begin{split}
\dfrac{d}{dt} v &= F(v) + \mu(I_h(u) - I_h(v)),\\
v(0)&=v_0.
\end{split}
\end{empheq}
This system, used in conjunction with \eqref{ODE}, is the AOT algorithm for data assimilation of system \eqref{ODE}.  In the case where the dynamical system \eqref{ODE} is the 2D Navier-Stokes equations, it was proven in \cite{Azouani_Olson_Titi_2014,Azouani_Titi_2014} that,  for any divergence-free initial data $v_0\in L^2$, 
$
 \|u(t) - v(t)\|_{L^2} \maps 0,
$
exponentially in time. 
In particular, even without knowing the initial data $u_0$, the solution $u$ can be approximately reconstructed for large times.  We emphasize that, as noted in \cite{Azouani_Olson_Titi_2014}, the initial data for \eqref{AOT} can be any $L^2$ function, even $v_0=0$. Thus, no information about the initial data is required to reconstruct the solution asymptotically in time.

The principal aim of this article is to develop a new class of nonlinear algorithms for data assimilation.  The main idea is to use a nonlinear modification of the AOT algorithm for data assimilation to try to drive the algorithm toward the true solution at a faster rate.  In particular, for a given, possibly nonlinear function $\mathcal{N}$, we consider a modification of \eqref{AOT} in the form:
\begin{empheq}[left=\empheqlbrace]{align}\label{AOT_NL}
\begin{split}
\dfrac{d}{dt} v &= F(v) + \mu \mathcal{N}(I_h(u) - I_h(v)),\\
v(0)&=v_0.
\end{split}
\end{empheq}
To begin, we first focus on the following form of the nonlinearity:
\begin{align}\label{NL_power}
\mathcal{N}(x) =\mathcal{N}_1(x) := x|x|^{-\gamma}, \quad x\neq0, \quad 0<\gamma<1,
\end{align}
with $\mathcal{N}_1(0)=0$.  
\begin{remark}\label{concave_remark}
 Note that by formally setting $\gamma=0$, one recovers the linear AOT algorithm \eqref{AOT}.  The main idea behind using such a nonlinearity is that, when $I_h(u)$ is close to $I_h(v)$, the solution $v$ is driven toward the true solution $u$ more strongly than in the linear AOT algorithm. In particular, for any $c>0$, if $x>0$ is small enough, then $\mathcal{N}_1(x)> cx$, so no matter how large $\mu$ is chosen in the linear AOT algorithm, the nonlinear method with $\mathcal{N}=\mathcal{N}_1$ will always penalize small errors more strongly.  
\end{remark}


As a preliminary test of the effectiveness of this approach, in this work we demonstrate the nonlinear data assimilation algorithm \eqref{AOT_NL} on a one-dimensional PDE; namely, the Kuramoto-Sivashinky equation (KSE), given in dimensionless units by:
\begin{empheq}[left=\empheqlbrace]{align}\label{KSE}
\begin{split}
   u_t + uu_x + \lambda u_{xx}+ u_{xxxx}&=0,
   \\
   u(x,0) &=u_0(x),
\end{split}
\end{empheq}
in a periodic domain $\Omega = [-L/2,L/2] = \nR/L\nZ$ of length $L$. Here, $\lambda>0$ is a dimensionless parameter.  For simplicity, we assume that the initial data is sufficiently smooth (made more precise below) and  mean-free, i.e., $\int_{-L/2}^{L/2}u_0(x)\,dx=0$, which implies $\int_{-L/2}^{L/2}u(x,t)\,dx=0$ for all $t\geq0$. This equation has many similarities with the 2D Navier-Stokes equations. It is globally well-posed; it has chaotic large-time behavior; and it has a finite-dimensional global attractor, making it an excellent candidate for studying large-time behavior.  It governs various physical phenomena, such as the evolution of flame-fronts, the flow of viscous fluids down inclined planes, and certain types of crystal growth (see, e.g., \cite{Kuramoto_Tsuzuki_1976,Sivashinsky_1977,Sivashinsky_1980_stoichiometry}).  Much of the theory of the 1D Kuramoto-Sivashinsky equation was developed in the periodic case in  \cite{Collet_Eckmann_Epstein_Stubbe_1993_Attractor, Constantin_Foias_Nicolaenko_Temam_1989_IM_Book, Goodman_1994,  Tadmor_1986, Temam_1997_IDDS,Tadmor_1986,Ilyashenko_1992} (see also \cite{Azouani_Olson_Titi_2014,Bronski_Gambill_2006,Cheskidov_Foias_2001,Foias_Kukavica_1995,Foias_Nicolaenko_Sell_1985,Foias_Nicolaenko_Sell_Temam_1988,Goldman_Josien_Otto_2015,Golovin_Davis_Nepomnyashchy_2000,Hyman_Nicolaenko_Zaleski_1986,Jardak_Navon_Zupanski_2010,Jolly_Kevrekidis_Titi_1990,Jolly_Rosa_Temam_2000,Kuramoto_Tsuzuki_1976,Lions_Perthame_Tadmor_1994,Liu_1991_KSE_Gevrey,Molinet_2000,Otto_2009,Sivashinsky_1977,Sivashinsky_1980_stoichiometry,Hyman_Nicolaenko_1986}).   For a discussion of other boundary conditions for \eqref{KSE}, see, e.g., \cite{Kuramoto_Tsuzuki_1976,Sivashinsky_1977, Robinson_2001, Pokhozhaev_2008,Larios_Titi_2015_BC_Blowup}.   
%
Discussions about the numerical simulations of the KSE, can be found in, e.g.,  \cite{Foias_Jolly_Kevrekidis_Titi_1991_Nonlinearity,Foias_Titi_1991_Nonlinearity,Foias_Jolly_Kevrekidis_Titi_1994_KSE,Lopez_1994_pseudospectral_KSE,Kalogirou_Keaveny_Papageorgiou_2015}.  Data assimilation in several different contexts for the 1D Kuramoto-Sivashinsky equation was investigated in \cite{Jardak_Navon_Zupanski_2010,Lunasin_Titi_2015}, who also recognized its potential as an excellent test-bed for data assimilation.

Using the nonlinear data assimilation algorithm \eqref{AOT_NL} in the setting of the Kuramoto-Sivashinsky equation, with the choice \eqref{NL_power} for the nonlinearity for some $0<\gamma<1$, we arrive at:
\begin{empheq}[left=\empheqlbrace]{align}\label{KSE_DA}
\begin{split}
   v_t + vv_x + v_{xx}+ v_{xxxx}&=\mu\,\,\text{sign}(I_h(u) - I_h(v))|I_h(u) - I_h(v)|^{1-\gamma},
   \\
   v(x,0) &=v_0(x).
\end{split}
\end{empheq}
We take $I_h$ to be orthogonal projection onto the first $c/h$ Fourier modes, for some constant $c$.  Other physically relevant choices of $I_h$, such as a nodal interpolation operator, have been considered in the case of the linear AOT algorithm (see, e.g., \cite{Gesho_Olson_Titi_2015}).  

\begin{remark}
In view of the ODE example $y'=y^{1/2}$, $y(0)=0$, which has multiple solutions, one might wonder about the well-posedness of equation \eqref{KSE_DA}.  While lack of well-posedness is a possibility, our simulations do not appear to be strongly affected by such a hindrance.  In any case, we believe the greatly reduced convergence time we observe makes the equations worth studying.  A similar remark can be made for an equation we examine in a later section with power larger than one, in view of the ODE example $z'=z^2$, $z(0)=1$, which develops a singularity in finite time.  A study of the well-posedness of \eqref{AOT_NL}, in various settings, will be the subject of a forthcoming work.
\end{remark}

\section{Preliminaries}\label{secPre}
\noindent

In this work, we compute norms of the error given by the difference between the data assimilation solution and the reference solution.  We focus on the $L^2$ and $H^1$ norms, defined by
\begin{align*}
 \|u\|_{L^2}^2 = \frac{1}{L}\int_{-L/2}^{L/2}|u(x)|^2\,dx,
 \qquad
 \|u\|_{H^1} = \|\nabla u\|_{L^2}.
\end{align*}
(Note that, by Poincar\'e's inequality, $\|u\|_{L^2}\leq C\|u\|_{H^1}$, which holds on any bounded domain, $\|u\|_{H^1}$ is indeed a norm.)

We briefly mention the scaling arguments used to justify the form \eqref{KSE}.  For $a>0$, $b>0$, $L>0$, consider an equation in the form
\begin{align}
 u_t + uu_x + au_{xx} + bu_{xxxx}=0,\qquad x\in[-L/2,L/2]
\end{align}
Choose time scale $T=L^4/b$, characteristic velocity $U=L/T=b/L^3$, and dimensionless number  $\lambda = aL^2/b$.  Write $u' = u/U$, $x'=(x+L/2)/L$, $t' = t/T$, where the prime denotes a dimensionless variable.  Then
\begin{align}\label{KSE_lambda}
 \frac{U}{T}u'_{t'} + \frac{U^2}{L}u'u'_{x'} + \frac{aU}{L^2}u'_{x'x'} + \frac{bU}{L^4}u'_{x'x'x'x'}=0,\qquad x'\in[0,1]
\end{align}
Multiply by $L^4/(bU)$.  The equation in dimensionless form then becomes 
\begin{align}
u'_{t'} + u'u'_{x'} + \lambda u'_{x'x'} + u'_{x'x'x'x'}=0,\qquad x'\in[0,1].
\end{align}
Thus, $\lambda$ acts as a parameter which influences the dynamics, in the same way that the Reynolds number influences dynamics in turbulent flows.

Another approach is to set $\ell = (b/a)^{1/2}$, $T=\ell^4/b=b/a^2$, and $U = \ell/T = a^{3/2}/b^{1/2}$.  Then define dimensionless quantities (denoted again by primes) $u' = u'/U$, $x'=(x+L/2)/\ell$, $t' = t/T$.  The equation now becomes
\begin{align}\label{KSE_lambda_1}
 \frac{U}{T}u'_{t'} + \frac{U^2}{\ell}u'u'_{x'} + \frac{aU}{\ell^2}u'_{x'x'} + \frac{bU}{\ell^4}u'_{x'x'x'x'}=0,\qquad x'\in[0,\tfrac{L}{\ell}]
\end{align}
Multiplying by $\ell^4/(bU)$ yields
\begin{align}
u'_{t'} + u'u'_{x'} + u'_{x'x'} + u'_{x'x'x'x'}=0,\qquad x'\in[0,\tfrac{L}{\ell}].
\end{align}
Thus, equation \eqref{KSE_lambda_1} is similar to equation \eqref{KSE_lambda}, with $\lambda=1$, except that the dynamics are influenced by the dimensionless parameter $L/\ell$.  In particular, the dynamics can be thought of as influenced by parameter $\lambda$ with $L$ fixed, or equivalently influenced by the length of the domain $L$ with $\lambda$ fixed, where $\lambda\sim (L/\ell)^2$.  In this work, for the sake of matching the initial data used in \cite{Kassam_Trefethen_2005}, we choose the domain to be $[-16\pi,16\pi]$, so $L=32\pi$ is fixed, and we let $\lambda$ be the parameter affecting the dynamics.

 \section{Computational Results}\label{secComputations}

In this section, we demonstrate some computational results for the nonlinear data assimilation algorithm given by \eqref{KSE_DA}.

\subsection{Numerical Methods}\label{subsecNumericalMethods}
It was observed in \cite{Gesho_Olson_Titi_2015} that no higher-order multi-stage Runge-Kutta-type method exists for solving \eqref{KSE_DA} due to the need to evaluate at fractional time steps, for which the data $I_h(u)$ is not available.  Therefore, we use a semi-implicit spectral method with Euler time stepping.  The linear terms are treated via a first-order exponential time differencing (ETD) method (see, e.g., \cite{Kassam_Trefethen_2005} for a detailed description of this method).  The nonlinear term is computed explicitly, and in the standard way, i.e.,  by computing the derivatives in spectral space, and products in physical space, respecting the usual 2/3's dealiasing rule.  We use $N=2^{13}=8192$ spatial grid points on the interval $[-16\pi,16\pi)$, so $ \Delta x= 32\pi/N \approx 0.0123$.  We use a fixed time-step respecting the advective CFL; in fact, we choose $\Delta t =1.2207\times 10^{-4}$. 
For simplicity, we choose $\mu = 1$, however, the results reported here are qualitatively similar for a wide range of $\mu$ values.  For example, when $\mu=10$, convergence times are shorter for all methods, but the error plots are qualitatively similar.  
In \cite{Kassam_Trefethen_2005}, the case $\lambda=1$ is examined.  However, to examine a slightly more chaotic setting, we take $\lambda=2$, which is still well-resolved with $N=8192$.  Our results are qualitatively similar for smaller values of $\lambda$.  

Here, we let $I_h$ be the projection onto the lowest $M=\lfloor 32\pi/h\rfloor$ Fourier modes.  In this work, we set $M=32$ (i.e, $h = \pi$); so only the lowest 32 modes of $u$ are passed to the assimilation equation via $I_h(u)$.  
One can consider a more general interpolation operator as well, such as nodal interpolation, but we focus on projection onto low Fourier modes.


To fix ideas, in this paper we mainly use the initial data used in \cite{Kassam_Trefethen_2005} to simulate \eqref{KSE}; namely
\begin{align}\label{KSE_ic}
 u_0(x) =  \cos(x/16)(1+\sin(x/16));
\end{align}
on the interval $[-16\pi,16\pi]$.  However, we also investigated several other choices of initial data.  In all cases, the results were qualitatively similar to the ones reported here.  We present one such test near the end of this paper.

Note that explicit treatment of the term $\mu(I_h(u) - I_h(v))$ imposes a constraint on the time step, namely $\Delta t<2/\mu$ (which follows from a standard stability analysis for Euler's method).  This is not a series restriction in this work, since we choose $\mu=1$.

\begin{figure}[htp]
\begin{center}
\subfigure[Time-averaged spectrum of the reference solution $u$ to the 1D-Kuramoto-Sivashinsky equation.]
{\includegraphics[width=.49\textwidth]{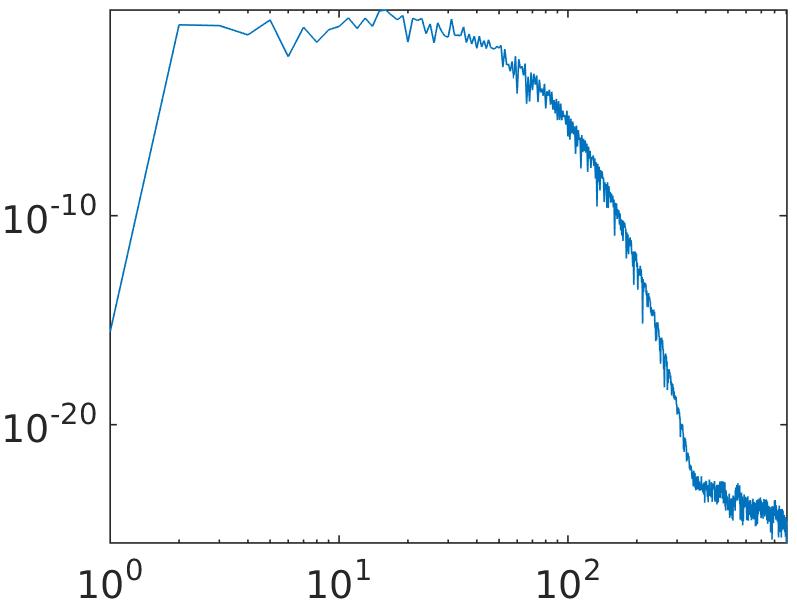}\label{FIG_spec}}
\hfill
\subfigure[Time-averaged spectrum of the data assimilation solution $v$ with nonlinear pure-power  ($\mathcal{N}_1$) algorithm.]
{\includegraphics[width=.49\textwidth]{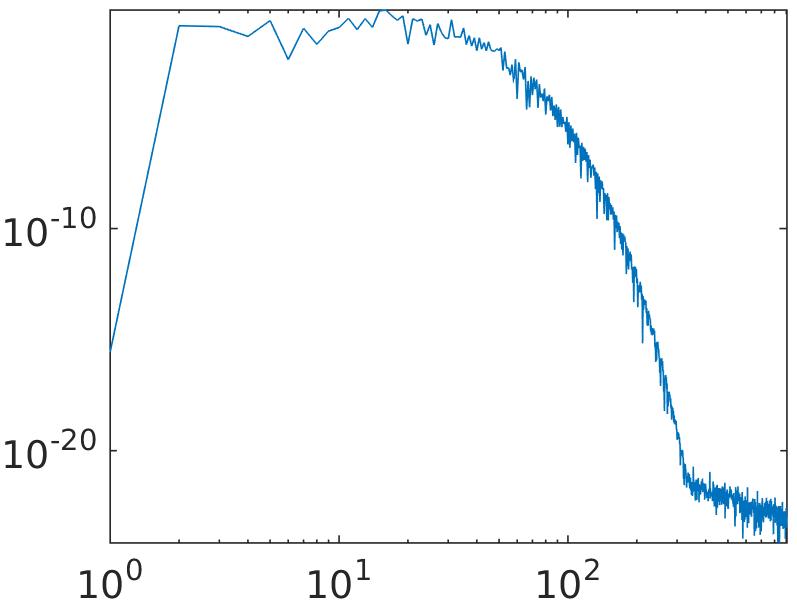}\label{FIG_spec_DA_v}}
\hfill
\subfigure[Time-averaged spectrum of the data assimilation solution $v$ with nonlinear hybrid ($\mathcal{N}_2$) algorithm.]
{\includegraphics[width=.49\textwidth]{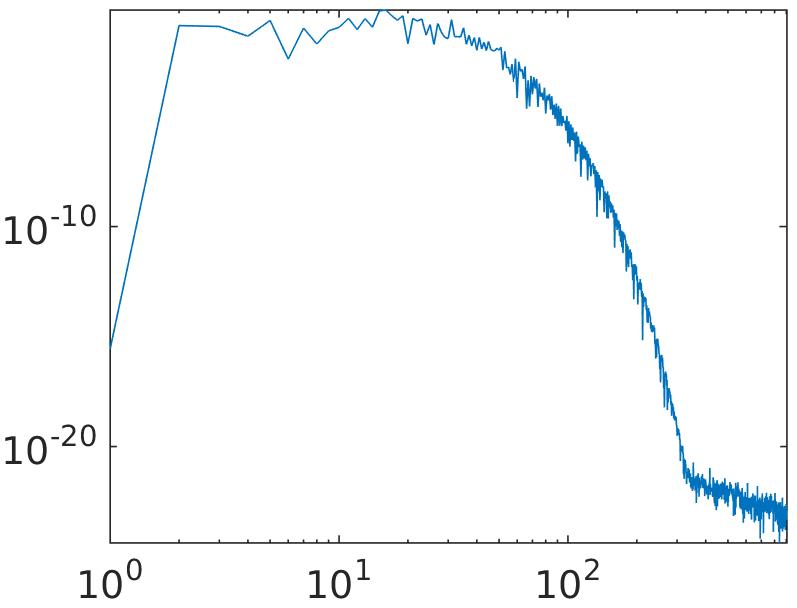}\label{FIG_spec_DAH_v}}
\hfill
\subfigure[Time-averaged spectrum of the data assimilation solution $v$ with nonlinear concave-convex ($\mathcal{N}_3$) algorithm.]
{\includegraphics[width=.49\textwidth]{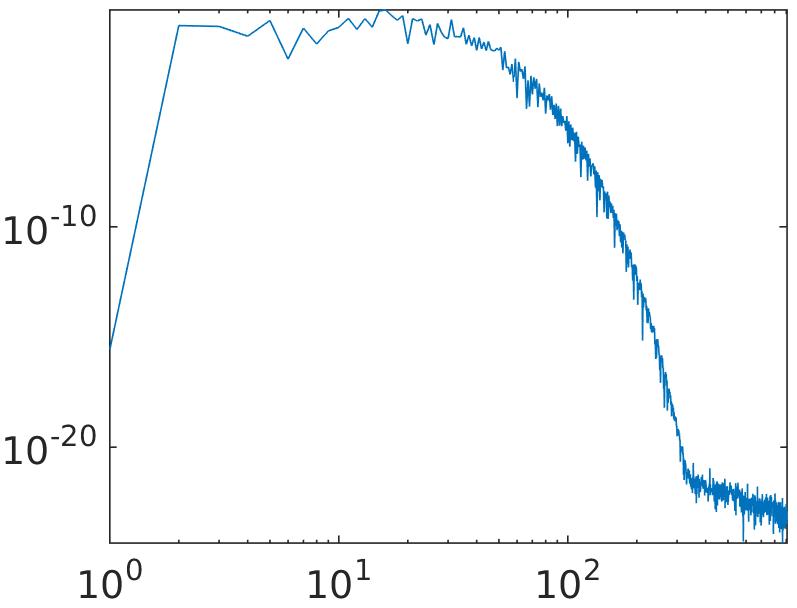}\label{FIG_spec_DANL_v}}
\caption{\label{figs:KSE_DA_Spectrum}\footnotesize 
 Log-log plots of the spectra for the above scenarios.  Plots are averaged over all time steps between times $t=20$ and $t=60$. ($\lambda=2$)} 
 \end{center} 
\end{figure}

All the simulations in the present work are well-resolved.  In Figure \eqref{figs:KSE_DA_Spectrum} we show plots of time-averaged spectra of all the PDEs simulated in the present work.  One can see that all relevant wave-modes are captured to within machine precision.  

\subsection{Simple Power Nonlinearity}
\begin{figure}[htp]
\begin{center}
\subfigure[A chaotic solution to the Kuramoto-Sivashinsky equation evolving in time.]
{\includegraphics[width=.32\textwidth]{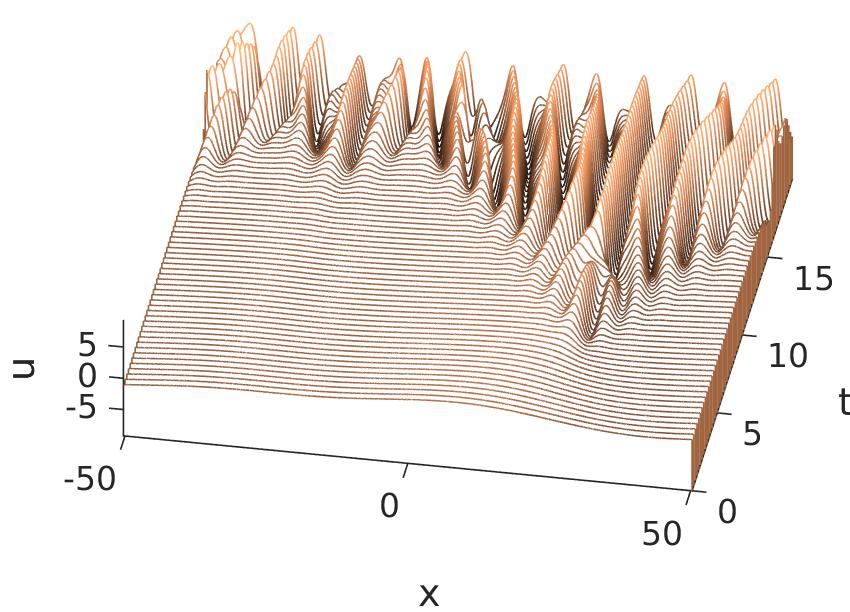}}
\hfill
\subfigure[Error in data assimilation solution using linear AOT algorithm ($\gamma=0$).]
{\includegraphics[width=.32\textwidth]{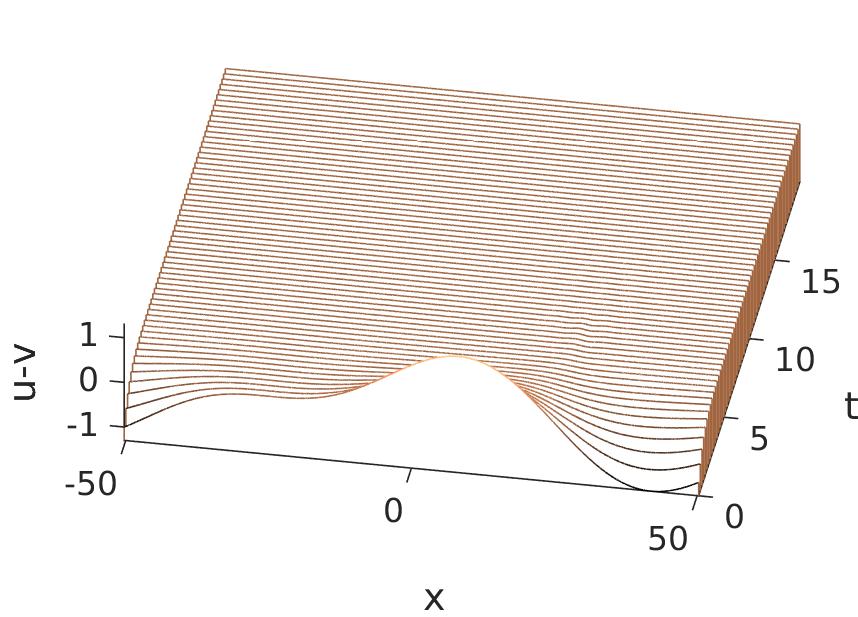}}
\hfill
\subfigure[Error in data assimilation solution using nonlinear algorithm \eqref{AOT_NL} ($\gamma=0.125$).]
{\includegraphics[width=.32\textwidth]{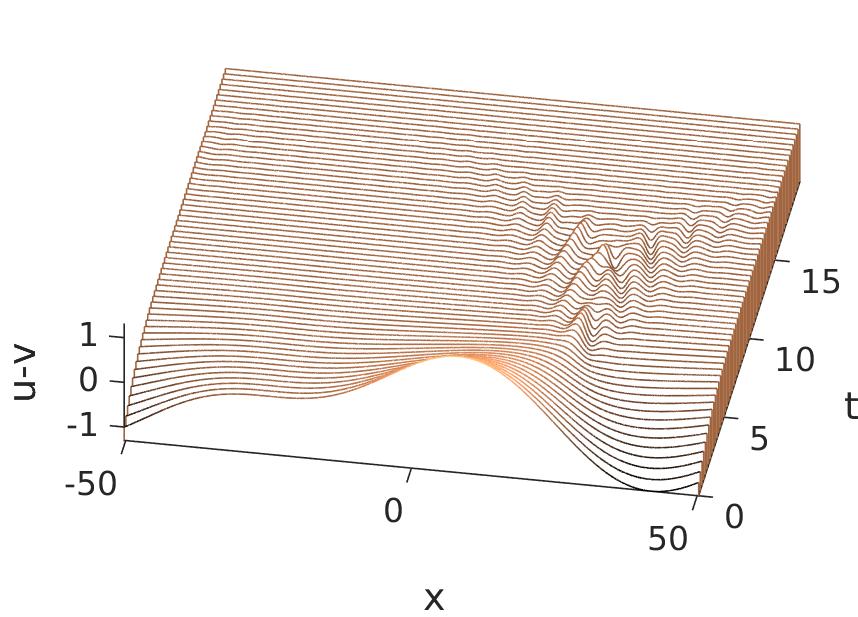}}
\caption{\label{figs:KSE_waterfall}\footnotesize 
 Data Assimilation for the Kuramoto-Sivashinsky equation ($\lambda=2$) using linear and nonlinear algorithms.  
 The difference rapidly decays to zero in time, and visually the errors look similar.  The assimilation equations were initialized with $v_0(x)=0$. $I_h$ is the orthogonal projection onto the lowest $32$ Fourier modes.  Similar results appear in tests of a wide variety of initial data, and for $0<\gamma<0.125$.}
\end{center}
\end{figure}

We compare the error in the nonlinear data assimilation algorithm \eqref{AOT_NL} with the error in the linear AOT algorithm. 
We first focus on nonlinearity given by a power according to \eqref{NL_power}; i.e., we consider equation \eqref{KSE_DA} together with equation \eqref{KSE}.  
In Figure \ref{figs:KSE_waterfall}(a), the solution to \eqref{KSE} (which we call the ``reference'' solution) evolves from the smooth, low-mode initial condition \eqref{KSE_ic} to a chaotic state after about time $t=20$.  In Figure \ref{figs:KSE_waterfall}(b), the difference between this solution and the AOT data assimilation solution is plotted.  It rapidly decays to zero in a short time.  

\begin{figure}[htp]
\begin{center}
\subfigure[Errors in $L^2$-norm vs. time.]
{\includegraphics[width=.49\textwidth]{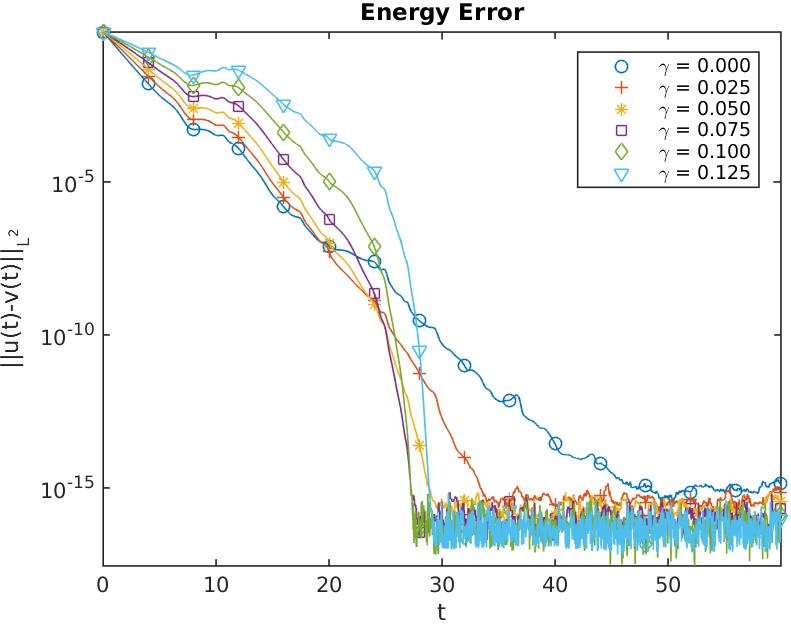}}
\hfill
\subfigure[Errors in $H^1$-norm vs. time.]
{\includegraphics[width=.49\textwidth]{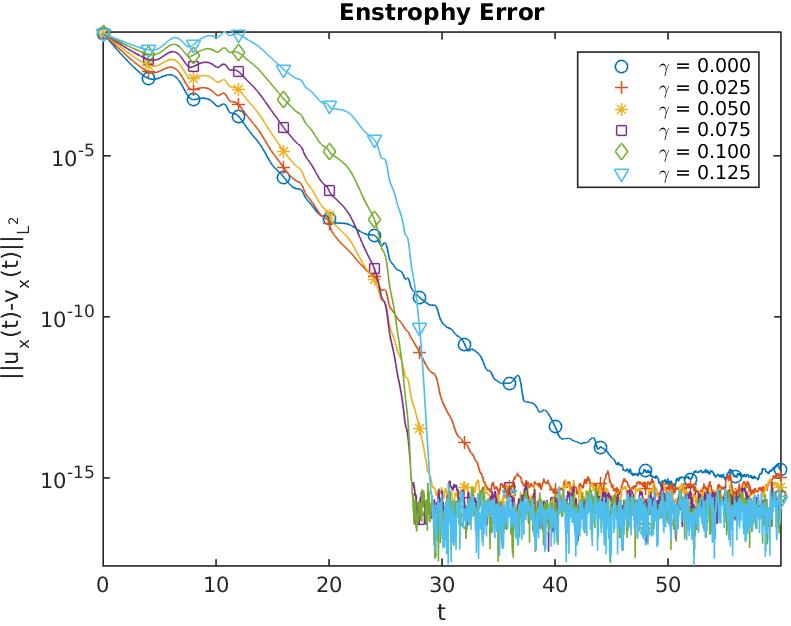}}
\caption{\label{figs:KSE_error}\footnotesize Error for the linear AOT ($\gamma=0$) solution and the nonlinear \eqref{AOT_NL} ($\gamma>0$) solution for various values of $\gamma$. Resolution 8192. (Log-linear scale.) }
\end{center}
\end{figure}

We observe in Figure \ref{figs:KSE_error} that errors in the linear AOT algorithm \eqref{AOT} and the nonlinear algorithm \eqref{AOT_NL} solutions both decay.  The error in the nonlinear algorithm has oscillations for roughly $5\lesssim t\lesssim 15$ which are not present in the error for the AOT algorithm.  However, by tracking norms of the difference of the solutions, one can see in Figure \ref{figs:KSE_error} that the nonlinear algorithm reaches machine precision significantly faster than the linear AOT algorithm, for a range of $\gamma$ values.  
When $\gamma>0.2$, our simulations appear to no longer converge (not shown here).  The error in the linear AOT algorithm (i.e.,  $\gamma=0$) reaches machine precision at roughly time $t\approx49.8$.  For $0<\gamma<0.2$, there seems to be an optimal choice in our simulations around $0.75\lesssim\gamma\lesssim 0.1$, reaching machine precision around $t\approx27.3$, a speedup factor of roughly $49.8/27.3\approx 1.8$.  Moreover, the shape of the curves with $\gamma>0$ indicate super-exponential convergence, indicated by the concave curve on the log-linear plot in Figure \ref{figs:KSE_error}, while for the linear AOT algorithm, the convergence is only exponential, indicated by the linear shape on the log-linear plot.  Currently, the super exponential convergence is only an observation in simulations.  An analytical investigation of the convergence rate will be the subject of a forthcoming work.  

\subsection{Hybrid Linear/Nonlinear Methods}\label{subsecHybrid}
In this subsection, we investigate a family of hybrid linear/nonlinear data assimilation algorithms. 
One can see from Figure \ref{figs:KSE_error}) in the previous subsection that, although the nonlinear methods converge to machine precision at earlier times than the linear method, the nonlinear method suffers from larger errors than the linear method for short times.  This motivates the possibility of using a hybrid linear/nonlinear method.  For example, one could look for an optimal time to switch between the models, say, perhaps around time $t\approx18\pm2$, according to Figure \ref{figs:KSE_error}), but this seems highly situationally dependent and difficult to implement in general.  Instead, the approach we consider here is to let $\mathcal{N}(x)$ be given by \eqref{NL_power} for $|x|\leq1$ but let it be linear for $|x|>1$.  The idea is that, when the error is small, deviations are strongly penalized, as in Remark \eqref{concave_remark}.  However, where the error is large, the linear AOT algorithm should give the greater penalization (i.e., $\mathcal{N}_1(x)<x$ when $x>1$).
Therefore, we consider algorithm \eqref{AOT_NL} with the following choice of nonlinearity, for some choice of $\gamma$, $0<\gamma<1$ (we take $\gamma=0.1$ in all following simulations).

\begin{align}\label{NL_hybrid}
\mathcal{N}(x) =\mathcal{N}_2(x) := 
\begin{cases}
 x,& |x|\geq 1,\\
 x|x|^{-\gamma}, & 0< |x|< 1,\\
 0,& x = 0.
\end{cases}
\end{align}
\begin{figure}[htp]
\begin{center}
\subfigure[Error in $L^2$-norm vs. time.]
{\includegraphics[width=.49\textwidth]{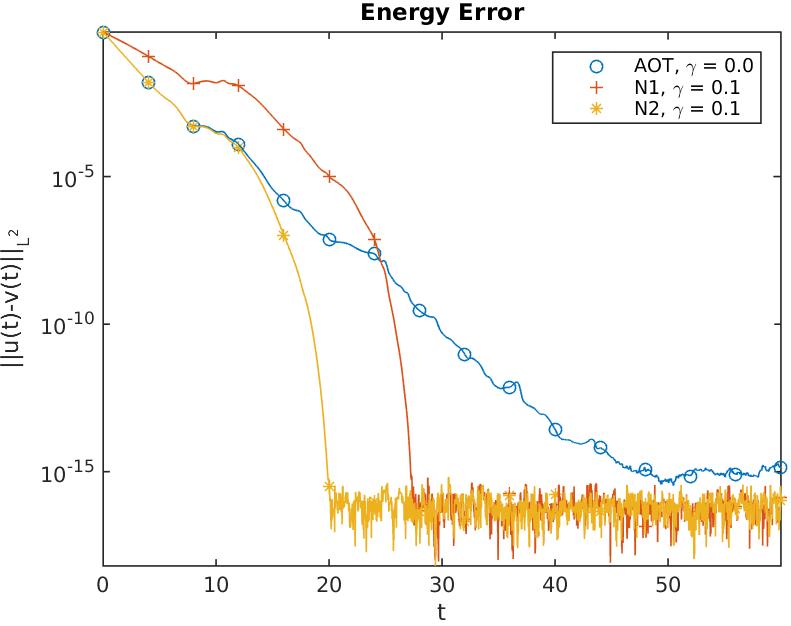}}
\hfill
\subfigure[Error in $H^1$-norm vs. time.]
{\includegraphics[width=.49\textwidth]{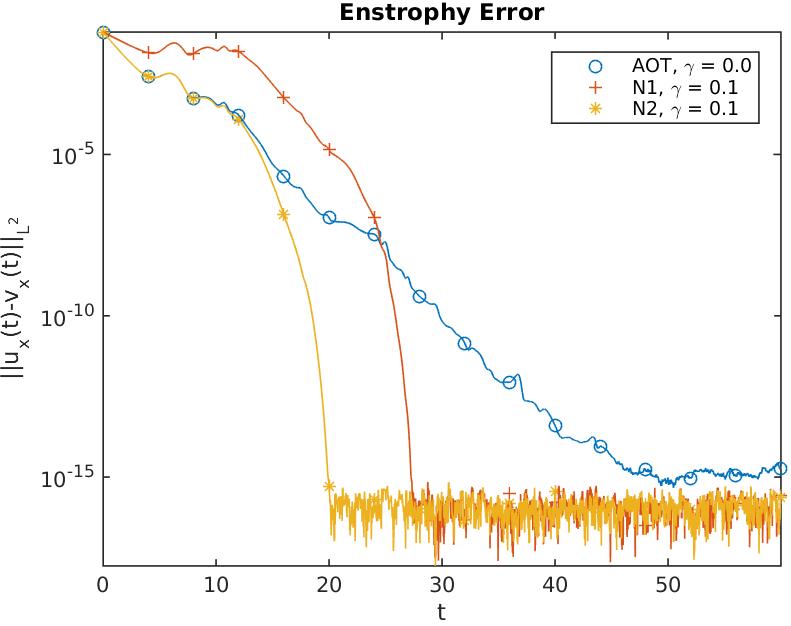}}
\caption{\label{figs:KSE_DAH_error}\footnotesize Error in linear ($\gamma=0$), nonlinear ($\gamma=0.1$), and hybrid ($\gamma=0.1$) algorithms ($\lambda=2$).  Resolution 8192.  (Log-linear scale.)}
\end{center}
\end{figure}

In Figure \ref{figs:KSE_DAH_error},we compare the linear algorithm \eqref{AOT} with the nonlinear algorithm \eqref{AOT_NL} with pure-power nonlinearity $\mathcal{N}_1$, given by \eqref{NL_power}, and also with hybrid  nonlinearity $\mathcal{N}_2$,  given by \eqref{NL_hybrid}.  The convergence to machine precision happens approximately at $t\approx49.8$ (for AOT), $t\approx27.3$ (for $\mathcal{N}_1$), and $t\approx20.0$ (for $\mathcal{N}_2$), respectively.  In addition, one can see that the hybrid algorithm remains close to the linear AOT algorithm for short times. Moreover, after a short time, the hybrid algorithm undergoes super-exponential convergence, converging faster than every algorithm analyzed so far.  The benefits of this splitting of the nonlinearity between $|x|>1$ and $|x|<1$ seem clear.  Moreover, this approach can be exploited further, which is the topic of the next subsection.

\subsection{Concave-Convex Nonlinearity}\label{subsecHybrid_NL}
\begin{figure}[htp]
\begin{center}
\subfigure[Error in $L^2$-norm vs. time.]
{\includegraphics[width=.49\textwidth]{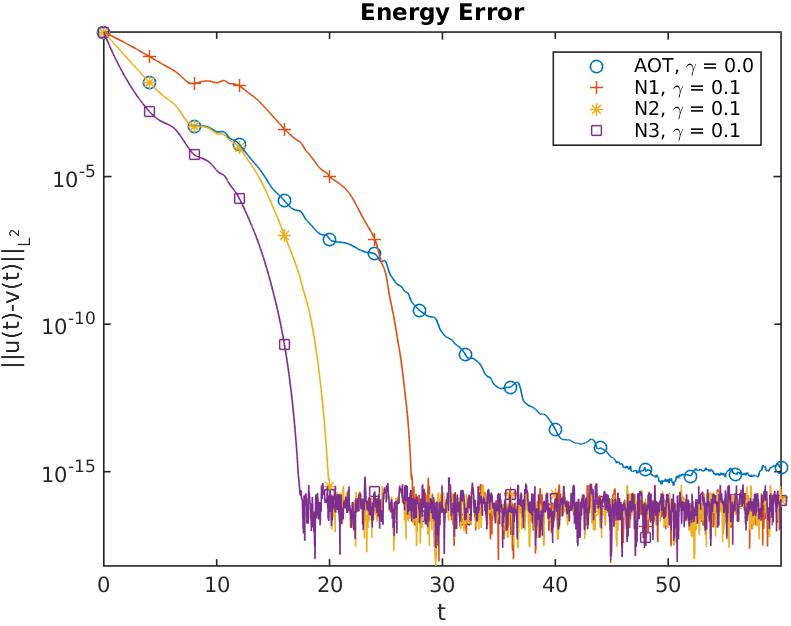}}
\hfill
\subfigure[Error in $H^1$-norm vs. time.]
{\includegraphics[width=.49\textwidth]{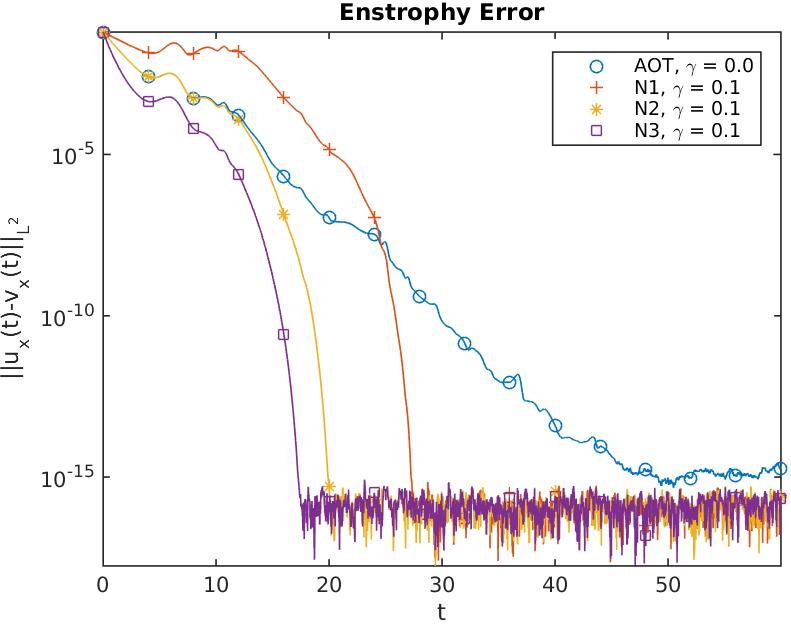}}
\caption{\label{FIG_err_L2_H1_all}\footnotesize Error in linear AOT  ($\gamma=0$) algorithm, and the non-linear algorithm with nonlinearities $\mathcal{N}_1$, $\mathcal{N}_2$, and $\mathcal{N}_3$ (each with $\gamma=0.1$).  ($\lambda=2$)  Resolution 8192. (Log-linear scale.)}
\end{center}
\end{figure}

Inspired by the success of the hybrid method, in this subsection, we further exploit the effect of the feedback control term $\mu\mathcal{N}(I_h(u)-I_h(v))$ by accentuating the nonlinearity for $|x|>1$.  We consider the following nonlinearity in conjunction with \eqref{AOT_NL} for the Kuramoto-Sivashinky equation. 
\begin{align}\label{NL_CC}
\mathcal{N}(x) =\mathcal{N}_3(x) := 
\begin{cases}
 x|x|^{\gamma},& |x|\geq 1,\\
 x|x|^{-\gamma}, & 0< |x|< 1,\\
 0,& x = 0.
\end{cases}
\end{align}
Note that this choice of $\mathcal{N}_3$ is concave for $|x|< 1$, and convex for $|x|\geq 1$.  The convexity for $|x|\geq 1$ serves to more strongly penalize large deviations from the reference solution.  In Figure \ref{FIG_err_L2_H1_all}, we see that \textit{at every positive time} this method has significantly smaller error than the linear AOT method, and the methods involving $\mathcal{N}_1$ and $\mathcal{N}_2$.  Convergence to machine precision happens at roughly $t\approx17.4$, a speedup factor of roughly $49.8/17.4\approx 2.8$ compared to the linear AOT algorithm.

\subsection{Comparison of All Methods}\label{subsecComparison}

\begin{figure}[htp]
\begin{center}
\subfigure[$t = 4$]
{\includegraphics[width=.49\textwidth]{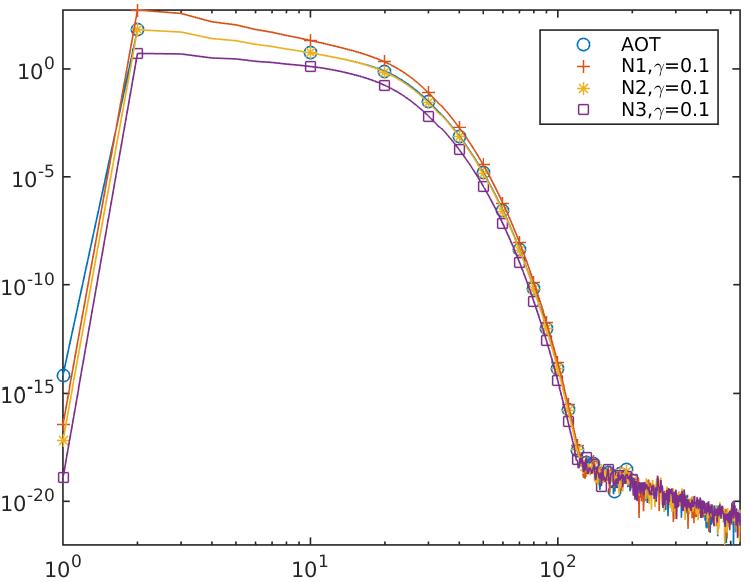}\label{FIG_err_spec4}}
\hfill
\subfigure[$t = 14$]
{\includegraphics[width=.49\textwidth]{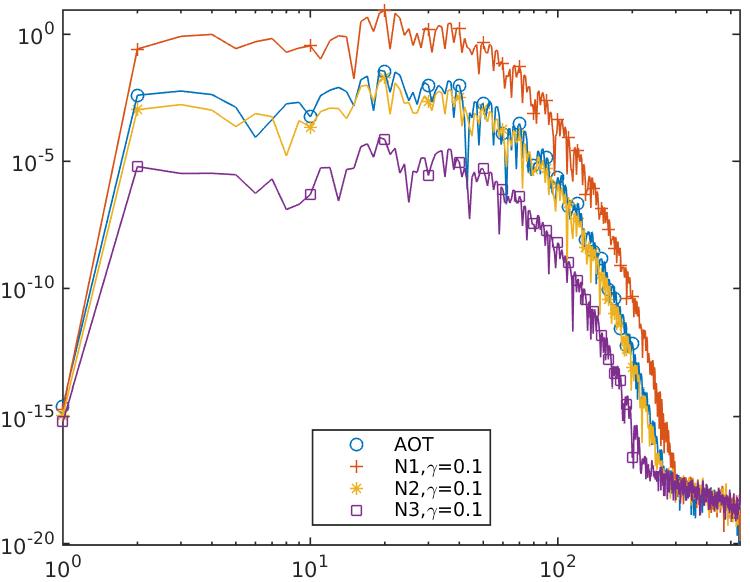}\label{FIG_err_spec14}}
\hfill
\subfigure[$t = 24$]
{\includegraphics[width=.49\textwidth]{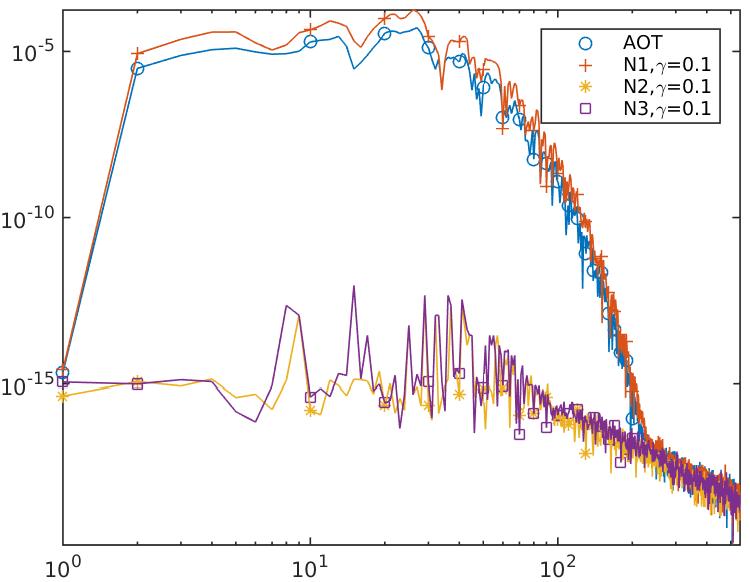}\label{FIG_err_spec24}}
\hfill
\subfigure[$t = 34$]
{\includegraphics[width=.49\textwidth]{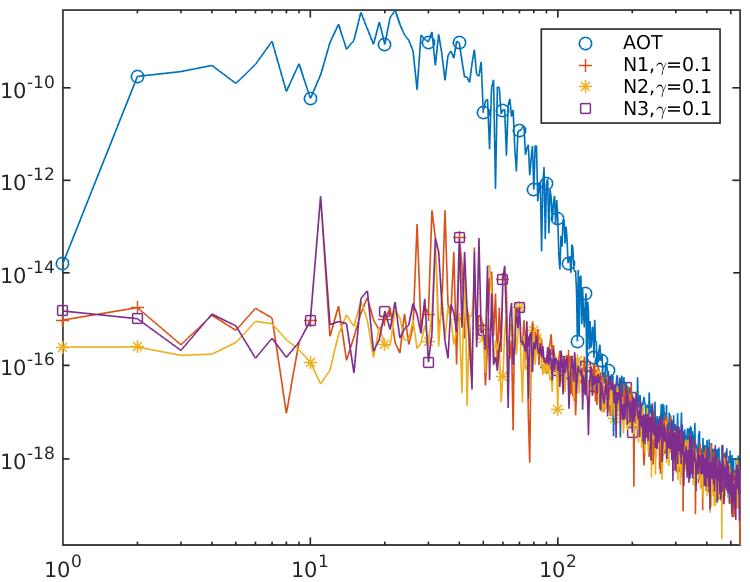}\label{FIG_err_spec34}}
\caption{\label{FIG_err_spec_all}\footnotesize 
Error in spectrum (mode amplitude vs. wave number) at different times for all methods.}
\end{center}
\end{figure}

Let us also consider the error at every Fourier mode. In Figure \ref{FIG_err_spec_all}, one can see these errors at various times.  We examine a time before the transition to fully-developed chaos ($t\approx 4$), at a time during the transition ($t\approx 14$), a time after the solution has settled down to an approximately statistically steady state  ($t\approx 24$), and a later time  ($t\approx 34$).   At each mode, and at each positive time, the error in the solution with nonlinearity \eqref{NL_CC} is the smallest. 

\begin{figure}[htp]
\begin{center}
\subfigure[Error in $L^2$-norm vs. time.]
{\includegraphics[width=.49\textwidth]{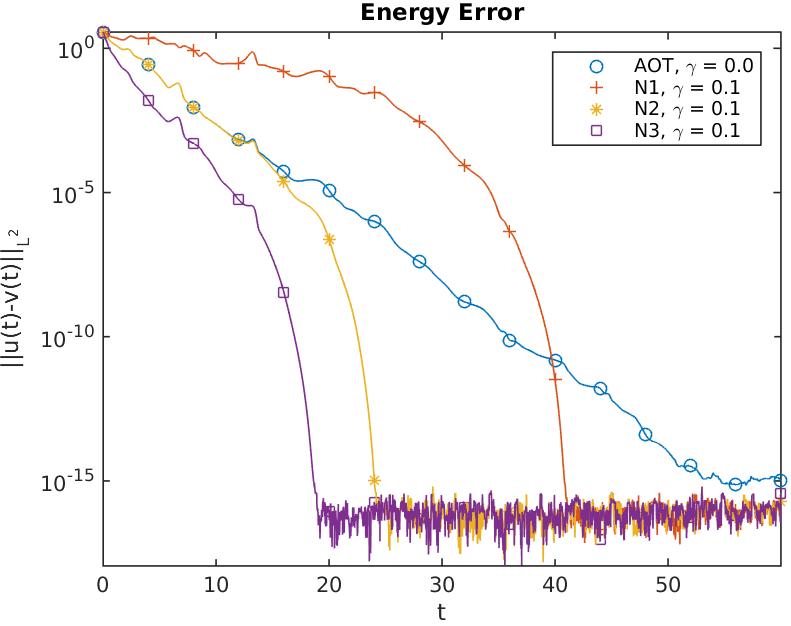}}
\hfill
\subfigure[Error in $H^1$-norm vs. time.]
{\includegraphics[width=.49\textwidth]{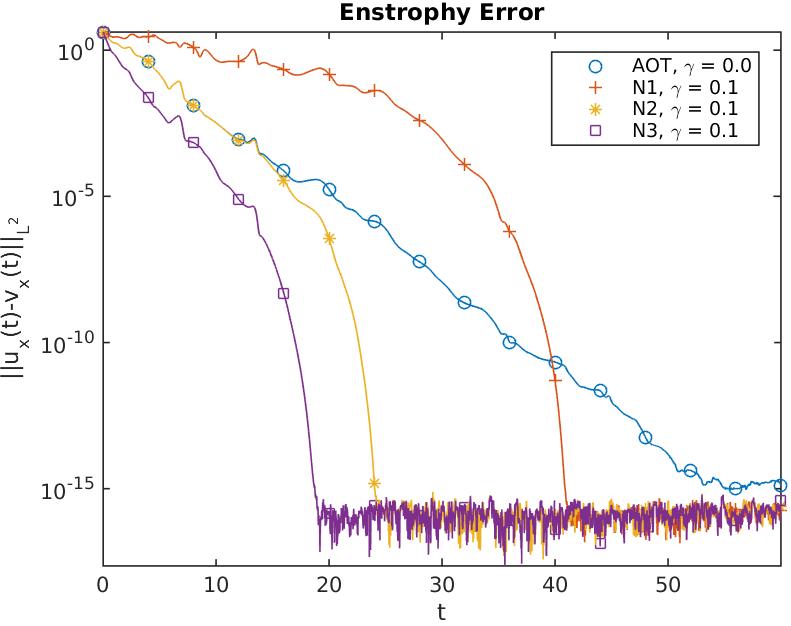}}
\caption{\label{FIG_err_chaos_L2_H1_all}\footnotesize Error in all algorithms with chaotic initialization for reference solution. ($\lambda=2$)  Resolution 8192. (Log-linear scale.)}
\end{center}
\end{figure}

Next, we point out that our results hold qualitatively with different choices of initial data for the reference equation \eqref{KSE}.  Therefore, we wait until the solution to \eqref{KSE} with initial data \eqref{KSE_ic} has reached an approximately statistically steady state (this happens roughly at $t\approx20$).  Then, we use this data to re-initialize the solution to \eqref{KSE} (in fact, we use the solution at $t=30$ to be well within the time interval of fully developed chaos).  We still initialize \eqref{KSE_DA} with $v_0\equiv0$. Norms of the errors are shown in Figure \ref{FIG_err_chaos_L2_H1_all}.  We observe that, although convergence time is increased for all methods, the qualitative observations discussed above still hold.

\section{Conclusions}\label{secConclusions}

Our results indicate that advantages might be gained by looking at nonlinear data assimilation.  We used the Kuramoto-Sivashinky equation as a proof-of-concept for this method; however, in a forthcoming work, we will extend the method to more challenging equations, including the Navier-Stokes equations of fluid flow.  Mathematical analysis of these methods will also be subject of future work.

We note that other choices of nonlinearity may very well be useful to consider.  Indeed, one may imagine a functional given by
\begin{align}
 \mathfrak{F}(\mathcal{N}) = t_*
\end{align}
where $t_*$ is the time of convergence to within a certain error tolerence, such as to within machine precision.  (One would need to show that admissible functions $\mathcal{N}$ are in some sense independent of the parameters and initial data, say, after some normalization.)  One could also consider a functional whose value at $\mathcal{N}$ is given by a particular norm of the error.  By minimizing such functionals, one might discover even better data-assimilation methods.

 
\bibliographystyle{abbrv}

\end{document}